\documentclass[oneside,11pt]{amsart}
\usepackage{amssymb,latexsym,amsmath,amsthm,enumitem,mathrsfs}
\usepackage[margin=1in]{geometry}
\usepackage{hyperref}
\usepackage{fancyhdr}
\usepackage{color}
\usepackage{stmaryrd}
\usepackage{mathrsfs}
\usepackage{lineno}

\newtheorem{letterthm}{Theorem}

\numberwithin{equation}{section}

\newtheorem{thm}{Theorem}[section]
\newtheorem*{thm*}{Theorem}
\newtheorem{prop}[thm]{Proposition}
\newtheorem{lemma}[thm]{Lemma}

\newtheorem{condition}[thm]{Condition}

\theoremstyle{remark}

\newcommand{\F}{\mathbb{F}}

\newcommand{\R}{\mathbb{R}}

\newcommand{\Z}{\mathbb{Z}}

\newcommand{\Cbf}{\mathbf{C}}
\newcommand{\Fscr}{\mathscr{F}}

\newcommand{\ep}{\varepsilon}

\newcommand{\con}{\equiv}

\newcommand{\ndiv}{\nmid}
\newcommand{\modd}[1]{\; ( \text{mod} \; #1)}
\newcommand{\bstack}[2]{#1 \atop #2}

\newcommand{\maps}{\rightarrow}
\newcommand{\intersect}{\cap}

\newcommand{\union}{\cup}

\newcommand{\al}{\alpha}
\newcommand{\be}{\beta}

\newcommand{\ga}{\gamma}
\newcommand{\del}{\delta}
\newcommand{\Del}{\Delta}

\newcommand{\Sig}{\Sigma}
\newcommand{\sig}{\sigma}

\newcommand{\Lam}{\Lambda}

\newcommand{\Acal}{\mathcal{A}}
\newcommand{\Pcal}{\mathcal{P}}

\newcommand{\Abf}{\mathbf{A}}
\newcommand{\Bbf}{\mathbf{B}}
\newcommand{\Lbf}{\mathbf{L}}
\newcommand{\Hbf}{\mathbf{H}}

\newcommand{\Qbf}{\mathbf{Q}}

\newcommand{\Ubf}{\mathbf{U}}

\newcommand{\Dbf}{\mathbf{D}}

\newcommand{\Nbf}{\mathbf{N}}
\newcommand{\Kbf}{\mathbf{K}}

\newcommand{\Xbf}{\mathbf{X}}

\newcommand{\zerobf}{\boldsymbol0}

\newcommand{\abf}{{\bf a}}

\newcommand{\kbf}{\mathbf{k}}

\newcommand{\mbf}{{\bf m}}
\newcommand{\nbf}{{\bf n}}

\newcommand{\tbf}{{\bf t}}

\newcommand{\xbf}{{\bf x}}

\newcommand{\beq}{\begin{equation}}
\newcommand{\eeq}{\end{equation}}

\makeatletter
\def\@tocline#1#2#3#4#5#6#7{\relax
  \ifnum #1>\c@tocdepth % then omit
  \else
    \par \addpenalty\@secpenalty\addvspace{#2}%
    \begingroup \hyphenpenalty\@M
    \@ifempty{#4}{%
      \@tempdima\csname r@tocindent\number#1\endcsname\relax
    }{%
      \@tempdima#4\relax
    }%
    \parindent\z@ \leftskip#3\relax \advance\leftskip\@tempdima\relax
    \rightskip\@pnumwidth plus4em \parfillskip-\@pnumwidth
    #5\leavevmode\hskip-\@tempdima
      \ifcase #1
       \or\or \hskip 1em \or \hskip 2em \else \hskip 3em \fi%
      #6\nobreak\relax
    \hfill\hbox to\@pnumwidth{\@tocpagenum{#7}}\par% <---- \dotfill -> \hfill
    \nobreak
    \endgroup
  \fi}
\makeatother

\newcommand{\Gbf}{\mathbf{G}}

\begin{document}

\title[]{Burgess bounds for short character sums evaluated at forms II: the mixed case}

\author[Pierce]{Lillian B. Pierce}
\address{Department of Mathematics, Duke University, 120 Science Drive, Durham NC 27708 USA}
\email{pierce@math.duke.edu}

 %Estimates on character sums

\begin{abstract}
 This work proves a Burgess bound for short mixed character sums in $n$ dimensions. The non-principal multiplicative character of prime conductor $q$ may be evaluated at any ``admissible'' form, and the additive character may be evaluated at any real-valued polynomial. The resulting   upper bound for the mixed character sum is nontrivial when the length of the sum  is at least $q^{\be}$ with $\be> 1/2 - 1/(2(n+1))$ in each coordinate. This work capitalizes on the recent stratification of multiplicative character sums due to Xu, and the resolution of the Vinogradov Mean Value Theorem in arbitrary dimensions.
\end{abstract}

\maketitle 

\section{Introduction}
Let $\chi$ be a non-principal multiplicative Dirichlet character modulo a prime $q$. Let $g \in \R[x_1,\ldots, x_n]$ be a polynomial of total degree $d \geq 1$, and let $F \in \Z[x_1,\ldots,x_n]$ be a form of degree $D \geq 1$. 
Define 
\[ S(F,g;\Nbf,\Hbf) = \sum_{\xbf \in (\Nbf,\Nbf+\Hbf]} e(g(\xbf))\chi(F(\xbf)),\]
where $\Nbf = (N_1,\ldots, N_n),$ $\Hbf = (H_1,\ldots, H_n)$ and $\xbf \in (\Nbf,\Nbf + \Hbf]$ denotes those tuples $\xbf \in \Z^n$ such that $x_i \in (N_i, N_i + H_i]$ for each $ 1\leq i \leq n$.
Such character sums are the building blocks of many methods in analytic number theory. 
The trivial bound is $|S(F,g;\Nbf,\Hbf)|  \leq H_1 \cdots H_n$, and   bounds that improve on this have many applications.  Conjecturally one could expect square-root cancellation to hold, for appropriate functions $F$ and $g$. In the particular case of short sums, namely those in which $H_i < q^{1/2}$, this remains out of reach, and a central goal is to provide any nontrivial upper bound, that is $|S(F,g;\Nbf,\Hbf)| = o(H_1 \cdots H_n)$, valid for general choices of $F,g$.

Historically the most fundamental case has been that of a one-dimensional multiplicative character sum, in which case Burgess's work set the gold standard, also establishing a long-standing subconvexity result for Dirichlet $L$-functions; see e.g.  \cite{Bur57,Bur63A}. 
(This subconvexity bound has only now been improved, in \cite{PetYou19x}.) Burgess's   method of proof has been resistant to substantial improvement, but recent work has begun to generalize the method to new settings. For a survey of Burgess bounds, in particular in the case of purely multiplicative sums, we refer to the overview given in \cite{PieXu19}.

In this paper we prove a  Burgess bound for mixed sums of the form $S(F,g;\Nbf,\Hbf)$, for the largest class of forms $F$ (acting nontrivially on all variables) for which one would anticipate a nontrivial bound  could be obtained. We formally define this set of ``admissible forms,'' before stating our main result. For this purpose, we recall that a polynomial $h$ is said to be $\Del$-th power-free over $\F_q$ if when $h$ is factored over $\F_q$ into irreducible pairwise non-associate factors $h_i$, each $h_i$ appears to a power strictly smaller than $\Del$.

  \begin{condition}[$(\Del,q)$-admissible]\label{cond} Fix a prime $q$ and an integer $\Del \geq 1$. A polynomial $f \in \F_q[x_1,\ldots, x_n]$ is $(\Del,q)$-admissible if the following holds. Upon writing $f = g^\Del h$ where $g,h \in \F_q[x_1,\ldots, x_n]$ and $h$ is $\Del$-th power-free over $\F_q$, then $h$  cannot be made independent of  a variable after a linear transformation, i.e.  there exists no  $A \in \mathrm{GL}_n(\F_q)$  such that $h(\xbf A) \in \F_q[x_2,\ldots, x_n]$. 
      \end{condition}

For any $\Del \geq 2$, if a form $F \in \Z[X_1,\ldots, X_n]$ satisfies $F=G^\Del H$ with $G,H \in \Z[x_1,\ldots, x_n]$ where $H$ is $\Del$-th power-free over $\Z$ and $H$ cannot be made  independent of a variable after a $GL_n(\Z)$ change of variables, then $F$ has $(\Del,q)$-admissible reduction modulo $q$ for all but finitely many primes $q$. An example of such a form is $x_1^D + \cdots +x_n^D$, and moreover, such forms are generic amongst the set of all forms in $\Z[x_1,\ldots,x_n]$ of degree $D$.   See  \cite[\S 3.1]{PieXu19} for further details on these facts.

Our main result is the following theorem.
\begin{thm}\label{thm_mixed}
Fix $n \geq 2$ and $d,D \geq 1$. Let $q$ be a fixed prime, and let $\chi$ be a non-principal Dirichlet character of  conductor $q$ and order $\Del$. 
Let $F \in \Z[x_1,\ldots, x_n]$ be a form of degree $D$ such that its reduction modulo $q$ is $(\Del,q)$-admissible. 
Let $g \in \R[x_1,\ldots, x_n]$ be a   polynomial of total degree $d \geq 1$.
Define for each integer $r \geq 1$,
\beq\label{theta_thm}
 \Theta = \Theta_{n, r} =\left\lfloor \frac{r-1}{n-1} \right\rfloor , \qquad M=M_{d, n} = d { n+d \choose n} \frac{n}{n+1}.
 \eeq
Let $\Hbf =(H,\ldots,H)$. 
Then for  every integer $r \geq 1$ such that $\Theta   = \Theta_{n, r}>  M$ and   $H<q^{\frac{1}{2}+\frac{1}{4 (\Theta -M)}}$,   
\beq\label{thm_ineq}
  |S(F,g;\Nbf,\Hbf) | \ll  H^{n-\frac{n+1}{2r}} q^{\frac{n(\Theta - M) +1 }{4r(\Theta-M)}}q^\ep, 
  \eeq
for every $\ep>0$;  the implied constant may depend on $n,d,D, \Del, r,\ep$ but is   independent of $g, F$.
\end{thm}

\noindent  \emph{Remark:} 
In dimension $n=1$, an inequality analogous to Theorem \ref{thm_mixed} with $\Theta = \Theta_{1,r}= r$ and $F(x)=x$ was  proved by Heath-Brown and the author \cite{HBP15}.    Upon setting $\Theta_{1,r}=r$ in each instance where $\Theta = \Theta_{n, r}$ appears in this paper, the method of the present paper also recovers this case, but we focus on $n \geq 2$.
At the time it was published, some results in  \cite{HBP15} were conditional on the  Main Conjecture in the Vinogradov Mean Value  Method, which now has been proved \cite{Woo16,BDG16}.

Theorem \ref{thm_mixed} is the first  Burgess bound for mixed sums in dimensions $n \geq 2$ in which $F$ is allowed to be any admissible form. In this sense it is a natural sequel to the work of the author with Xu \cite{PieXu19}, which introduced the $(\Del,q)$-admissible class of forms, in the setting of purely multiplicative sums (namely the case $S(F,0;\Nbf,\Hbf)$). Theorem \ref{thm_mixed} is   of comparable strength to the purely multiplicative case considered in \cite{PieXu19}.
Precisely, define 
\beq\label{beta_dfn}
 \be_n = \frac{1}{2} - \frac{1}{2(n+1)}.
 \eeq
Theorem \ref{thm_mixed} provides a nontrivial bound of the form $|S(F,g;\Nbf,\Hbf)| \ll H^n q^{-\del}$ when $H > q^{\be_n + \kappa}$ for some sufficiently small $\kappa >0$, and the savings is of the strength 
\[ \del \approx \frac{(n+1)^2}{4(n-1)} \kappa^2 \]
as $\kappa \maps 0$; see \S \ref{sec_quant} for details.
In particular, note that this savings is independent of the degree $D$ of the form $F$ and the degree $d$ of the polynomial $g$; this is achieved by an application of the sharp upper bound in the multi-dimensional Vinogradov Mean Value Theorem, due to \cite{PPW13} in many cases and \cite{GuoZha18} in complete generality.

Earlier work on two special types of mixed sums in dimensions $n \geq 2$ appeared in two recent papers.
 In the  special case   $F(\xbf) = x_1 \cdots x_n$, the author  proved nontrivial bounds for $|S(F, g;\Nbf,\Hbf)|$ as long as $H_i > q^{1/4 + \kappa}$ for some small $\kappa >0$ \cite{Pie16}. See also the preprint of Kerr \cite{Ker14x} in the case that $F(\xbf) = \prod_{i=1}^n L_i(\xbf)$ is the product of $n$ linear forms $L_i$ that are linearly independent over $\mathbb{F}_q$. 
In each of these special settings, additional structure allowed the argument to achieve the Burgess threshold $q^{1/4 + \kappa}$ for any $\kappa>0$, in any dimension.

 \subsection{Method of proof}
The proof of Theorem \ref{thm_mixed} capitalizes upon recent foundational work of two kinds:
\begin{enumerate} 
\item Xu's stratification of multiplicative character sums   \cite{Xu18};
\item the resolution of the Main Conjecture in the setting of the Vinogradov Mean Value Method. In the ``one-dimensional''  setting this is due to Wooley \cite{Woo16} in the cubic case and to Bourgain, Demeter and Guth \cite{BDG16} for all higher-degree cases (see also Wooley \cite{Woo19}). In  the higher-dimensional setting, this is due to Guo and Zhang \cite{GuoZha18}, and in certain regimes the earlier work of Parsell, Prendiville, and Wooley \cite{PPW13}. 
\end{enumerate}

In particular, the proof of Theorem \ref{thm_mixed}  applies a sharp upper bound for the number of solutions to the (multi-dimensional) Vinogradov system 
\beq\label{Vin_sys_generic}
 \xbf_1^\be + \cdots + \xbf_r^\be = \xbf_{r+1}^\be + \cdots + \xbf_{2r}^{\be} , \qquad 1 \leq |\be| \leq d
 \eeq
with $\xbf_j \in \Z^n$ and  $1 \leq  x_{j,i} \leq X$ for $1 \leq j \leq 2r$, $1 \leq i \leq n$; here $\be = (\be_1,\ldots, \be_n)$ is a multi-index with $|\be| = \be_1 + \cdots + \be_n$.
 In fact our work also applies to more general translation-dilation invariant systems (see Theorem \ref{thm_mixed_T}), 
 and as our method naturally uses the properties of such systems, we introduce the relevant terminology in the following section.

\section{Introducing the associated Vinogradov system}

To prove Theorem \ref{thm_mixed} for a fixed (admissible) choice of $F \in \Z[x_1,\ldots, x_n]$ and $g \in \R[x_1,\ldots, x_n]$, our primary object of focus will be 
\[ T(F,\Gbf;\Nbf,\Hbf) = \sup_{g \in \Fscr_0(\Gbf)} \sup_{\Kbf \leq \Hbf} \left| \sum_{\xbf  \in (\Nbf,\Nbf+\Kbf]} e(g(\xbf)) \chi(F(\xbf)) \right|,\]
in which $\Fscr_0(\Gbf)$ is a certain set of polynomials including the fixed polynomial $g$ of our choice.
Our main bound for $T(F,\Gbf;\Nbf,\Hbf)$ involves counting the number of solutions to a system of Diophantine equations, which we now introduce precisely.

To prove Theorem \ref{thm_mixed}, we  may take $\Gbf$ to be the set of all non-constant monomials in $n$ variables of total degree at most $d$, that is 
  \beq\label{G_choice}
  \Gbf   = \{ \xbf^\be \in \Z[x_1,\ldots, x_n], \be = (\be_1,\ldots, \be_n) \in \Z_{\geq 0}^n, 1 \leq |\be | \leq d\},
  \eeq
  in which $|\be| = \be_1 + \cdots + \be_n$. 
 (Momentarily, we will also consider other sets of monomials.) Given any set $\Gbf$ of monomials, we define $\Fscr_0(\Gbf)$ to be the set of all real-variable polynomials that are linear combinations of the elements in $\Gbf \union \{1\}$ (that is, including constant terms). 
  We will call the set $\Gbf $ defined in (\ref{G_choice}) the \emph{standard system} of monomials  in $n$ variables of degree at most $d$. With this choice of $\Gbf$, given any polynomial $g \in \R[x_1,\ldots, x_n]$ of degree $d$, we can embed it in $\Fscr_0(\Gbf)$. In particular,
  $|S(F,g;\Nbf,\Hbf)|  \leq T(F,\Gbf;\Nbf,\Hbf) $.
  
  The main outcome of the Burgess argument we develop  is an upper bound for $T(F,\Gbf; \Nbf,\Hbf)$ in terms of  (i) a complete multiplicative character sum and (ii) a complete additive character sum. We apply Xu's stratification \cite{Xu18} to bound the complete multiplicative character sum.  
  We evaluate the complete additive character sum precisely, and then dominate the outcome by the number of integral solutions to the (multi-dimensional) Vinogradov system of Diophantine equations (\ref{Vin_sys_generic}) associated to the system $\Gbf$ given in (\ref{G_choice}). (This is sometimes also called a Parsell-Vinogradov system of equations when $n \geq 2$.)
  
 Let $J_r(\Gbf,X)$ denote the number of integral solutions to the system (\ref{Vin_sys_generic}) with $1 \leq  x_{j,i} \leq X$ for $1 \leq j \leq 2r$, $1 \leq i \leq n$. 
We also let $M=M(\Gbf)$ denote the sum of the total degrees of all multi-indices in $\{ \be \in \Z^n_{\geq 0}: 1 \leq |\be| \leq d\}$.

Our main result in the context of Theorem \ref{thm_mixed} is as follows. For any $n \geq 2$ and $r \geq 1$, define $\Theta = \Theta_{n, r}$ as in Theorem \ref{thm_mixed}. For all integers $r \geq n$ and any $\Hbf = (H,\ldots, H)$ and $P \leq H$ with $HP<q$ and $ P \leq H  q^{-1/2\Theta}$,  
\[
 T(F,\Gbf; \Nbf,\Hbf) \ll  (H/P)^{M/2r}  H^{-n/2r}P^{n-1/2r}q^{n/4r}( \log q)^{n+1} \\
 	   \{  J_r(\Gbf,2H/P)^{1/2r}  +  q^{1/4r} (H/P)^{n-\Theta/2r}\}.
\]
Theorem \ref{thm_mixed} then follows    from an appropriate bound for $J_r(\Gbf,X)$ provided by the multi-dimensional Vinogradov Mean Value Theorem, and an optimal choice for $P$ in terms of $H,q$.

\subsection{Remark on more general systems $\Gbf$}
 Without any additional difficulty, our main arguments can replace the standard system $\Gbf$ specified in (\ref{G_choice})   by any reduced monomial translation-dilation invariant system. This terminology was introduced in   \cite{PPW13}, and we briefly recall the definitions.
A given collection  $\Gbf = \{g_1, \ldots, g_R\}$  of $R$ non-constant monomials in $\Z[x_1,\ldots, x_n]$ is said to be translation-dilation invariant if there exist polynomials  $c_{m,\ell} \in \Z[\xi_1,\ldots, \xi_n]$ 
  for $1 \leq m \leq R, 0 \leq \ell \leq m$ with $c_{m,m}=1$ for $1 \leq m \leq R$ and such that for any $\xi \in \Z^n$, 
  \[ g_m(\xbf + \xi) = c_{m,0}(\xi) + \sum_{\ell=1}^m c_{m,\ell}(\xi) g_\ell(\xbf), \quad 1 \leq m \leq R.\]
   (See \cite[Eqn (2.3)]{PPW13} for an explanation of why such systems are called translation-dilation invariant.)
The system $\Gbf$ is said to be reduced if  the set $\{g_1, \ldots, g_R\}$ is linearly independent over $\R$.
  To avoid degenerate cases, we will only work with systems $\Gbf$ that include all variables nontrivially, and in particular include linear monomials in each variable.
   
For either the standard system  (\ref{G_choice}) or for any reduced monomial translation-dilation invariant system $\Gbf$, the  following quantities  will arise in our proof. Given $\Gbf$ as above, we say it has dimension $n$ and rank $R$.
 Let $\Lambda(\Gbf)$ be the associated set of multi-indices, so that $\Gbf=\{ \xbf^\be : \be \in \Lambda (\Gbf)\}$. 
 We set the degree  $d(\Gbf) = \max \{ |\be| : \be \in \Lambda(\Gbf)\}$ to be the highest total degree appearing in a monomial in $\Gbf$.
The rank is $R(\Gbf) = |\Lambda(\Gbf)|$ and we define the weight $M(\Gbf)$ (or homogeneous dimension) by 
\beq\label{rank_dfn}  M(\Gbf) = \sum_{\be \in \Lambda(\Gbf)} |\be| .
\eeq
For the standard system $\Gbf$ in (\ref{G_choice}) of monomials in $n$ variables of total degree at most $d$, 
  \beq\label{MRga}
  R= R(\Gbf)= { n+d \choose n} -1, \qquad  M =M(\Gbf) =  d { n+d \choose n} \frac{n}{n+1} .
  \eeq
We define the associated Vinogradov system of $R(\Gbf)$ equations in $2r$ variables by
   \beq\label{Vin_sys}
 \xbf_1^\be + \cdots + \xbf_r^\be = \xbf_{r+1}^\be + \cdots + \xbf_{2r}^{\be}, \qquad \be \in \Lambda (\Gbf). \eeq
We let 
  $J_r(\Gbf,X)$ denote the number of integral solutions to the system (\ref{Vin_sys}) with $1 \leq  x_{j,i} \leq X$ for $1 \leq j \leq 2r$, $1 \leq i \leq n$. 
 In full generality, our methods prove that $T(F,\Gbf;\Nbf,\Hbf)$ can be controlled by the number of solutions $J_r(\Gbf,X)$.
    \begin{prop}\label{prop_T}
Let $n \geq 2$. Let $q$ be a fixed prime and let $\chi$ be a non-principal Dirichlet character of conductor $q$ and order $\Delta$. Let $F \in \Z[x_1,\ldots, x_n]$ be a form of degree $D$, with $(\Del,q)$-admissible reduction modulo $q$. Let $\Gbf $ be a reduced monomial translation-dilation invariant system (containing linear monomials in each variable) with weight $M(\Gbf)$. For each $r \geq 1$ define $\Theta = \Theta_{n, r} = \lfloor \frac{r-1}{n-1} \rfloor$. For all integers $r \geq n$ and any $\Hbf = (H,\ldots, H)$ and $P \leq H$ with $HP<q$ and $ P \leq H  q^{-1/2\Theta}$,  
\[ T(F,\Gbf; \Nbf,\Hbf) \ll  (H/P)^{M(\Gbf)/2r}  H^{-n/2r}P^{n-1/2r}q^{n/4r}( \log q)^{n+1}  
   \{  J_r(\Gbf,2H/P)^{1/2r}  +   q^{1/4r}  (H/P)^{n-\Theta/2r}\}.
\]
 \end{prop}
 Thus   for any translation-dilation invariant system $\Gbf$ for which a suitable bound is known for $J_r(\Gbf,X)$, we can deduce a Burgess bound for $|S(F, g; \Nbf,\Hbf)|$, for any polynomial $g$ in the span of $\Gbf$. 

\noindent  \emph{Remark:} Given $g$, one could optimize the choice of $\Gbf$ as in \cite{Pie16}, but we do not pursue this here.

\subsection{Key results for Vinogradov Mean Value Theorems in multi-dimensional settings}
Once we have proved Proposition \ref{prop_T}, it is clear that the key remaining step to prove Theorem \ref{thm_mixed} is to bound $J_r(\Gbf,X)$.  In dimensions $n \geq 2$,  Parsell, Prendiville, and Wooley \cite{PPW13} proved that 
 for any reduced translation-dilation invariant system $\Gbf$, for all $r > R(\Gbf) (d(\Gbf)+1)$,   the sharp upper bound for $J_r(\Gbf,X)$ holds. For this range of $r$, the sharp upper bound is $J_r(\Gbf,X)\ll_{r,n , d, \ep} X^{2nr-M(\Gbf) + \ep}$.
Recently, Guo and Zhang \cite{GuoZha18} have proved 
  the sharp upper bound for $J_r(\Gbf,X)$ for the standard system   (\ref{G_choice}), 
  for all $n \geq 2, d \geq 1$ and all $r \geq 1$; the exact form of the sharp upper bound depends on the size of $r$. 
  This completely resolves the Main Conjecture for the multi-dimensional Vinogradov Mean Value Method, for the standard system (\ref{G_choice}). (See also the earlier proof of the sharp upper bound for $n=2, d=2$ in \cite{BouDem16} and $n=2,d=3$ in \cite{BDG17}.)     
  More generally, Guo and Zorin-Kranich \cite{GuoZor18}  have now proved sharp upper bounds for $J_r(\Gbf,X)$ for all $r$, for any system $\Gbf$ that is an   Arkhipov-Chubarikov-Karatsuba system. We provide more details on these bounds for $J_r(\Gbf,X)$, and their implications for Burgess bounds, in an appendix in \S \ref{sec_gen}.

\section{Initiating the Burgess argument}
 
For the remainder of the paper, we assume that $n \geq 2$ and that a prime $q$ has been fixed; we then fix a non-principal multiplicative character $\chi$ of conductor $q$ and order $\Del$.  We fix a degree $D$ and assume $F \in \Z[x_1,\ldots,x_n]$ is a form of degree $D$ such that its reduction modulo $q$ is $(\Del,q)$-admissible. 
We let $\Gbf$ be a reduced monomial translation-dilation invariant system with all corresponding notation as defined above; to avoid degenerate situations, we assume that $\Gbf$ contains linear monomials in each of the $n$ variables. 
In particular, for Theorem \ref{thm_mixed}, we can take $\Gbf$ as in (\ref{G_choice}). 
We then define
\[ T(F,\Gbf;\Nbf,\Hbf) = \sup_{g \in \Fscr_0(\Gbf)} \sup_{\Kbf \leq \Hbf} \left| \sum_{\xbf  \in (\Nbf,\Nbf+\Kbf]} e(g(\xbf)) \chi(F(\xbf)) \right|.\]

The construction  $T(F,\Gbf;\Nbf,\Hbf)$ (which also appeared in \cite{HBP15}) has several advantageous properties in comparison to a sum $S(F, g; \Nbf,\Hbf)$ with a fixed polynomial $g$. First,   $T(F,\Gbf;\Nbf,\Hbf)$ is periodic under any shift of $\Nbf$ by multiples of $q$, and thus we will assume from now on that $0 \leq N_i < q$ for $i =1,\ldots, n$.  Second, the fact that $T(F,\Gbf;\Nbf,\Hbf)$ includes a supremum over polynomials in $\Fscr_0(\Gbf)$ will allow us to replace a supremum over ranges of summation by a supremum over linear phases via Fourier inversion (Lemma \ref{lemma_sum}), which in turn  is subsumed in the supremum over polynomials in $\Fscr_0(\Gbf)$. Finally, and most crucially, the supremum over $g \in \Fscr_0(\Gbf)$ will allow us to run the Burgess argument including the factor $e(g(\xbf))$, as we now demonstrate.

To begin the Burgess argument, we suppose $\Hbf = (H, \ldots, H)$ is fixed with $H<q$. We consider any $\Kbf \leq \Hbf$, by which we mean $K_i \leq H$ for $i=1,\ldots, n$. For any tuple $\Kbf$ we denote $\|\Kbf\| = K_1\cdots K_n$. For a parameter $P$ assumed to satisfy $1 \leq P \leq H  $ we define the set $\Pcal$ of auxiliary primes by 
\[ \Pcal = \{ P < p \leq 2P : p \ndiv q\},\]
so that $|\Pcal| \gg P/\log P$.  
We write each $\xbf \in (\Nbf,\Nbf+ \Kbf]$ according to its residue class modulo $p$, as 
\[ \xbf = \abf q + \mbf p,\]
where $\abf = (a_1,\ldots, a_n)$ with $0 \leq a_i < p$ and $\mbf  \in (\Nbf^{\abf,p}, \Nbf^{\abf,p} + \Kbf^{p}]$, with the definitions 
\[\Nbf^{\abf,p} = \Nbf/p - \abf q/p, \qquad \Kbf^{p} = \Kbf/p. \]
Then by applying the periodicity and multiplicativity of $\chi$ and the homogeneity of $F$, for any $g \in \Fscr_0(\Gbf)$,
\begin{align*}
 \sum_{\xbf  \in (\Nbf,\Nbf+\Kbf]} e(g(\xbf)) \chi(F(\xbf))
	& = \sum_{\bstack{\abf}{ 0 \leq a_i < p}} \sum_{\mbf  \in (\Nbf^{\abf,p}, \Nbf^{\abf,p} + \Kbf^{p}]} 
	 	e(g(\abf q + \mbf p)) \chi(F(\abf q + \mbf p))\\
		& =  \chi(p^D)\sum_{\bstack{\abf}{ 0 \leq a_i < p}}  \sum_{\mbf  \in (\Nbf^{\abf,p}, \Nbf^{\abf,p} + \Kbf^{p}]} 
	 	e(g(\abf q + \mbf p)) \chi(F( \mbf)).
	 \end{align*}
In particular, we note that $\Kbf^{p} \leq \Kbf/p \leq \Hbf/p \leq \Hbf/P$.
Consequently, after taking absolute values and taking the supremum over $g \in \Fscr_0(G)$ and $\Kbf \leq \Hbf$, we have
\[ T(F,\Gbf;\Nbf,\Hbf) \leq \sum_{\bstack{\abf}{ 0 \leq a_i < p}} T(F,\Gbf;\Nbf^{\abf,p},\Hbf/P).\]
Finally, we average this inequality over all $p \in \Pcal$, so that
	\beq\label{T_ap}
	 T(F,\Gbf;\Nbf,\Hbf) \leq |\Pcal|^{-1} \sum_{p \in \Pcal}  \sum_{\bstack{\abf}{ 0 \leq a_i < p}}  T(F,\Gbf;\Nbf^{\abf,p},\Hbf/P).
	 \eeq
\emph{Remark:} In  \cite{Pie16}, we   restricted to the special case $F(x_1,\ldots,x_n) = x_1 \cdots x_n$,  and  we could freely average over a distinct set of primes in each coordinate. Due to averaging over a larger set, we could recover a nontrivial bound for $H_i$ as small as $q^{1/4+\kappa}$ for $\kappa>0$. In our present setting, we can see from the argument above that to exploit the homogeneity of $F$ we must use the same prime $p$ for each coordinate, leading to a smaller set to average over. Nevertheless, many of the arguments of \cite{Pie16} may be adapted, and thus we will be efficient in our presentation.

Our next step is to introduce further averaging so that we may free the starting points $\Nbf^{\abf,p}$ from the dependence on $\abf,p$, enabling us to later interchange the order of summation, and apply H\"older's inequality.  
\begin{lemma}
Fix  $\Ubf \in \R^n$ and $\Lbf \in \R^n_{\geq 1}$. For any $\Kbf \leq \Lbf$, 
\[ T(F,\Gbf;\Ubf,\Kbf) \leq 2^{2n} \| \Lbf\|^{-1} \sum_{\Ubf - \Lbf < \mbf \leq \Ubf} T(F,\Gbf;\mbf, 2\Lbf).\]
\end{lemma}
This follows verbatim from the inclusion-exclusion proof given in  \cite[Lemma 3.1]{Pie16}, with each instance of $\chi_1(x_1) \cdots \chi_n(x_n)$ replaced by 
$\chi(F(\xbf))$, so we do not repeat the proof here.  (See also \cite[Lemma 5.1]{PieXu19} for more details on the inclusion-exclusion.) 

We apply the lemma to  (\ref{T_ap}) with $\Lbf = \Hbf/P$ (recalling $H/P \geq 1$) to obtain
 \beq\label{T_ap'}
   T(F,\Gbf;\Nbf,\Hbf) \ll \| \Hbf/P\|^{-1} |\Pcal|^{-1} \sum_{p \in \Pcal} \sum_{\bstack{\abf}{0 \leq a_i < p}}  
 	\sum_{\Nbf^{\abf,p} - \Hbf/P< \mbf \leq \Nbf^{\abf,p}} T(F,\Gbf;\mbf,2\Hbf/P).\eeq
Now for each $\mbf$ we define 
\[ \Acal(\mbf) = \{p \in \Pcal, \abf, 0 \leq a_i < p : \Nbf^{\abf,p} -\Hbf/P < \mbf \leq \Nbf^{\abf,p} \}.\]
By Lemma 5.2 of \cite{PieXu19}, $\Acal(\mbf)$ vanishes unless $|m_i| \leq 2q$ for each $i$, and moreover 
as long as 
\beq\label{HP_assp}
HP < q,
\eeq
which we henceforward assume, then
\[ \sum_{\mbf} \Acal(\mbf) \ll \sum_{\mbf} \Acal(\mbf)^2 \ll P \| \Hbf \| .\]
Applying H\"older's inequality twice to (\ref{T_ap'}) then shows that 
\begin{multline*}
  T(F,\Gbf;\Nbf,\Hbf) \\
     \ll  \| \Hbf/P\|^{-1} |\Pcal|^{-1} (\sum \Acal(\mbf))^{1 -1/r} ( \sum \Acal(\mbf)^2)^{1/2r} \left( \sum_{\mbf, |m_i| < 2q} T(F,\Gbf;\mbf, 2\Hbf/P)^{2r} \right)^{1/2r}.
  \end{multline*}
After simplification (recalling the periodicity of $T(F,\Gbf;\mbf,\Hbf)$ under shifts of $\mbf$ modulo $q$), we see that
\beq\label{T_ap''}
  T(F,\Gbf;\Nbf,\Hbf) \ll  \| \Hbf \|^{-1/2r}P^{n-1/2r} (\log P) \left( \sum_{\mbf \modd{q}} T(F,\Gbf;\mbf, 2\Hbf/P)^{2r} \right)^{1/2r}.\eeq

\subsection{Strategy to remove the suprema}
We required the definition of $T(F,\Gbf;\Nbf,\Hbf)$ to include two suprema in order to complete the various averaging arguments in the opening steps of the Burgess method, described above. Now we work to remove these suprema, in order to reveal a complete character sum over $\mbf \modd{q}$. 

We define 
\[ T_0(F,\Gbf;\mbf,\Kbf) =   \sup_{g \in \Fscr_0(\Gbf)}  \left| \sum_{\xbf  \in (\mbf,\mbf+\Kbf]} e(g(\xbf)) \chi(F(\xbf)) \right|
	=   \sup_{g \in \Fscr_0(\Gbf)}  \left| \sum_{\xbf  \in (\zerobf,\Kbf]} e(g(\xbf)) \chi(F(\xbf+ \mbf)) \right|
	.\]
Next, we suppose that we have indexed a finite set of polynomials $\theta_\al \in \R[x_1,\ldots, x_n]$ according to a finite set of indices $\al$, and for each such polynomial we define
\[ T_1(F,\theta_\al;\mbf,\Kbf) = \left|\sum_{\zerobf < \xbf \leq \Kbf} e(\theta_{\al}(\xbf)) \chi(F(\xbf + \mbf)) \right| .\]
\noindent  \emph{Remark:} Note that this is equal to $S(F, g; \mbf,\Kbf)$ with $g (\cdot) = \theta_{\al}(\cdot - \mbf)$, but for technical reasons it is easier to work with the notation $T_1$, which builds the shift by $\mbf$ into the argument of the multiplicative character.

To remove the suprema in the expression $T$, we will pass from expressions involving $T$ to expressions involving $T_0$, and then to expressions involving $T_1$. Then we will be ready to evaluate the contribution of the additive character sum exactly, and to apply Xu's stratification  to bound the contribution of the multiplicative character sum.

\subsection{Approximations of the additive character contribution}
We  first pass from $T$ to $T_0$ inside (\ref{T_ap''}), by applying Lemma 3.3 of \cite{Pie16} (an $n$-dimensional version of \cite[Lemma 2]{BomIwa86}), which we recall here:
\begin{lemma}\label{lemma_sum}
Let $a(\nbf)$ be a sequence of complex numbers indexed by integral tuples $\nbf$ supported on the set $\nbf \in (\Abf,\Abf+\Bbf]\subset \Z^n.$ Let $I = (\Cbf,\Cbf + \Dbf]$ be any product of intervals with $I \subseteq (\Abf, \Abf+\Bbf]$. Then 
\[ \sum_{\nbf \in I} a(\nbf) \ll ( \prod_{i=1}^n \log (B_i +2)) \sup_{\theta \in \R^n} \left| \sum_{\nbf \in (\Abf, \Abf + \Bbf]} a(\nbf) e(\theta \cdot \nbf) \right|.\]
\end{lemma}
This lemma shows that for any $\mbf$,
\[  T(F,\Gbf;\mbf,2\Hbf/P) \ll (\prod_{i=1}^n \log (2H_i/P+2)) T_0(F,\Gbf;\mbf,2\Hbf/P) \ll (\log q)^nT_0(F,\Gbf;\mbf,2\Hbf/P),\]
since $2H_i/P < 2q$. 
Note that here we use the fact that $\Gbf$ contains linear monomials in each variable, so that the supremum over $\theta \in \R^n$ is subsumed in the supremum over $g \in \Fscr_0(\Gbf)$. In the setting of Theorem \ref{thm_mixed}, we are using the hypothesis  that the degree $d$ of the polynomial $g$ is at least 1.

 Applying this in (\ref{T_ap''}) proves 
\beq\label{T_0}
 T(F,\Gbf; \Nbf,\Hbf) \ll  \| \Hbf \|^{-1/2r}P^{n-1/2r}( \log q)^{n+1} \left( \sum_{\mbf \modd{q}} T_0(F,\Gbf;\mbf, 2\Hbf/P)^{2r} \right)^{1/2r}.
 \eeq

In order to pass from $T_0$ to expressions involving $T_1$, we must fix a  set of representative polynomials $\theta_\al$ (indexed by $\al$) with the following property: for each $\mbf$, one of the representative polynomials $\theta_\al $ (depending on $\mbf$) has the property that  $T_1(F,\theta_\al;\mbf,2\Hbf/P)$ is sufficiently close in value to $T_0(F,\Gbf;\mbf,2\Hbf/P)$.  

We define these representative polynomials as in \cite[\S 4]{Pie16}, according to a fixed integer $Q \geq 1$ (to be chosen later), and $\Qbf = (Q,\ldots, Q)$. 

 We let $\Lambda_0(\Gbf) = \Lambda(\Gbf) \union \{(0,\ldots, 0)\}$ so $|\Lambda_0(\Gbf)| = R+1 = R(\Gbf)+1$. Since within the phase of an exponential sum, the coefficients of a polynomial $g$ are regarded modulo 1, $\Fscr_0(\Gbf)$ is represented by $[0,1]^{R+1}$, and upon ordering the $R+1$ multi-indices $\be \in \Lambda_0(\Gbf)$ in a fixed manner as $\be^{(0)} = \mathbf{0}, \ldots, \be^{(R)}$ once and for all, we partition the $\be^{(j)}$-th unit interval $[0,1]$ in this product $[0,1]^{R+1}$ into $\Qbf^{\be^{(j)}}  = Q^{|\be^{(j)}|}  $ sub-intervals of length $Q^{-|\be^{(j)}|}$.

  We recall from \cite[\S 4]{Pie16} the following facts  about this decomposition. 
  Recall the weight $M=M(\Gbf)$.
The decomposition partitions $[0,1]^{R+1}$ into  $Q^{M}$ boxes, which we call $B_\al$, according to indices $\al$ that we order once and for all. For each such box $B_\al$ we assign its distinguished vertex $\theta_\al$ to be the vertex with the least value in each coordinate, which is of the form
\beq\label{theta_spec}
 \theta_\al = (\theta_{\al,\be^{(0)}}, \ldots, \theta_{\al, \be^{(R)}}) = (c_{\be^{(0)}}Q^{-|\be^{(0)}|} , \ldots, c_{\be^{(R)}}Q^{-|\be^{(R)}|}),
 \eeq
where $c_{\be^{(0)}}=0$ and for each $j=1, \ldots, R$, $c_{\be^{(j)}}$ is an integer with $0 \leq c_{\be^{(j)}} \leq Q^{|\be^{(j)}|} - 1$.
 
 \noindent  \emph{Remark:} The outcome here is simplified relative to \cite{Pie16} since in our setting all the coordinates of $\Qbf = (Q,\ldots,Q)$ are the same. In particular we do not require the notion of the ``density''  of the system $\Gbf$, introduced with the notation $\ga(\Gbf)$ in \cite{Pie16}; this was the sum of the exponents in $\Lam(\Gbf)$. A reader comparing the notation here to \cite{Pie16} will observe that any term of the form $\Qbf ^\ga$ in the previous work can be written here as $Q^M$ with $M=M(\Gbf)$ the weight of the system $\Gbf$.

Now for any   point $\theta \in [0,1]^{R+1}$, we define an associated real-valued polynomial in $\R[X_1,\ldots,X_n]$ by
\beq\label{theta_poly_dfn}
 \theta(\Xbf) := \sum_{\be \in \Lambda_0 (\Gbf) } \theta_\be \Xbf^\be.
 \eeq
 In particular, for each box $B_\al$ with distinguished vertex $\theta_\al \in [0,1]^{R+1}$, we define the associated polynomial $\theta_\al(\Xbf)$. 
 Finally, we define  
\[ S_F(\kbf) :=  \sum_{\al} \sum_{\mbf \modd{q}} T_1(F,\theta_\al;\mbf, \kbf)^{2r}.
\]
Here the sum over $\al$ denotes the sum over the finitely many indices $\al$ in the decomposition.
 The following lemma records an upper bound for $T(F,\Gbf; \Nbf,\Hbf)$ in terms of $S_F(\kbf)$, according to this decomposition. 
 
    \begin{lemma}\label{lemma_TS}
Fix any $Q \geq 2H/P$. Let $[0,1]^{R+1} = \union_\al B_\al$ be partitioned as described above, according to indices $\al$.  
Also assume $HP < q$.
Then 
 \beq\label{TS}
 T(F,\Gbf; \Nbf,\Hbf) \ll  \| \Hbf \|^{-1/2r}P^{n-1/2r}( \log q)^{n+1} \left(  \sup_{\kbf \leq 2\Hbf/P} S_F(\kbf) \right)^{1/2r}.
 \eeq
 
 \end{lemma}
The proof follows that of \cite[Lemma 4.3]{Pie16} verbatim, upon replacing each appearance of $\chi_1(x_1) \cdots \chi_n(x_n)$ by $\chi(F(\xbf))$. Thus we do not repeat the proof in full detail, but highlight the most important steps. Fix $\mbf$ and consider the corresponding term $T_0(F,\Gbf; \mbf, 2\Hbf/P)^{2r}$ on the right-hand side of (\ref{T_0}).  Since the coefficients of polynomials $g \in \Fscr_0(\Gbf)$ are regarded modulo 1,   by compactness the supremum over $g \in \Fscr_0(\Gbf)$  in $T_0(F,\Gbf;\mbf , 2\Hbf/P)$ occurs for a particular polynomial, say $\tilde{g}$ (depending on $\mbf$).  Then 
$
T_0(F,\Gbf;\mbf , 2\Hbf/P) = T_1(F,\tilde{g};\mbf, 2\Hbf/P).
$

The partition of $[0,1]^{R+1}$ constructed according to the parameter $Q$ contains a box $B_\al$ with  index $\al$ (depending on $\mbf$) with the following property:  for each multi-index $\be^{(j)}$ with $j=0,\ldots, R,$ if $\theta_{\al,\be^{(j)}}$ denotes the corresponding coordinate of the distinguished vertex in $B_\al$, and $\tilde{g}_{\be^{(j)}}$ denotes  the coefficient  of $\Xbf^{\be^{(j)}}$ in $\tilde{g}(\Xbf)$, then 
\[ | \theta_{\al,\be^{(j)}} - \tilde{g}_{\be^{(j)}}| \leq Q^{-|\be^{(j)}|}, \qquad  j=0, \ldots, R.\]
Under the assumption that $Q \geq 2H/P$, 
 partial summation shows that replacing the polynomial $\tilde{g}(\Xbf)$  in $T_1(F,\tilde{g};\mbf,2\Hbf/P)$ by the  polynomial $\theta_\al(\Xbf)$ corresponding to this box $B_\al$ makes a sufficiently small error.  To state this precisely, we define the notation that for any subset $J \subseteq \{1,\ldots,n\}$   with cardinality $0 \leq |J| \leq n$ and complement $^cJ = \{1 , \ldots, n\} \setminus J$, 
 \[ T_1^{(^cJ),(J)} (F,\theta_\al; \mbf, \kbf_{(^cJ)}, \tbf_{(J)}) = \left| \sum_{\bstack{0 < x_j  \leq k_j}{j \in ^cJ}} \sum_{\bstack{0 < x_j \leq t_j}{j \in J}} e(\theta_\al(\xbf)) \chi(F(\xbf + \mbf)) \right|.  \]
 Then partial summation (applied as in \cite[Lemma 4.1]{Pie16}) shows that as long as $Q \geq 2H/P$, 
 \[ T_1(F, \tilde{g}; \mbf, 2\Hbf/P) \ll \sum_{J \subseteq \{1,\ldots, n\}} (2H/P)^{-|J|} 
 	\int \cdots \int_{(0,2H/P]^{|J|}} T_1^{(^cJ),(J)} (F,\theta_\al; \mbf, (2\Hbf/P)_{(^cJ)}, \tbf_{(J)}) d \tbf_{(J)}.\]
 A repeated application of H\"older's inequality then shows that 
\begin{multline*}
T_1(F, \tilde{g}; \mbf, 2\Hbf/P)^{2r}  \\
\ll \sum_{J \subseteq \{1,\ldots, n\}} (2H/P)^{-|J|} 
 	\int \cdots \int_{(0,2H/P]^{|J|}} T_1^{(^cJ),(J)} (F,\theta_\al; \mbf, (2\Hbf/P)_{(^cJ)}, \tbf_{(J)})^{2r} d \tbf_{(J)}.
	\end{multline*}
Now, since we do not know which index $\al$ was chosen to approximate $\tilde{g}(\Xbf)$ by $\theta_\al(\Xbf)$,   we replace the right-hand side by the sum of this expression over all $\al$; by positivity, this only enlarges the right-hand side. We conclude that 
\begin{multline*}
 T_1(F, \tilde{g}; \mbf, 2\Hbf/P)^{2r} \\
  \ll \sum_{J \subseteq \{1,\ldots, n\}} (2H/P)^{-|J|} 
 	\int \cdots \int_{(0,2H/P]^{|J|}} \sum_{\al} T_1^{(^cJ),(J)} (F,\theta_\al; \mbf, (2\Hbf/P)_{(^cJ)}, \tbf_{(J)})^{2r} d \tbf_{(J)}.
	\end{multline*}
This statement now holds uniformly in $\mbf$, and we can sum it over all $\mbf \modd{q}$. 
Note that by positivity, for any $\tbf_{(J)} \in (0,2H/P]^{|J|},$ 
\[ \sum_\al \sum_{\mbf \modd{q}} T_1^{(^cJ),(J)} (F,\theta_\al; \mbf, (2\Hbf/P)_{(^cJ)}, \tbf_{(J)})^{2r}
    \leq  \sup_{\kbf \leq 2\Hbf/P}   \sum_\al \sum_{\mbf \modd{q}}  T_1  (F,\theta_\al; \mbf, \kbf)^{2r} .\]
In conclusion, we have proved that 
 \[ \sum_{\mbf \modd{q}}  T_1(F, \tilde{g}; \mbf, 2\Hbf/P)^{2r} \ll  \sup_{\kbf \leq 2\Hbf/P}   \sum_\al \sum_{\mbf \modd{q}}  T_1  (F,\theta_\al; \mbf, \kbf)^{2r}  ,\]
 and this suffices to complete the proof of Lemma \ref{lemma_TS}.

\emph{Remark:}
While the introduction of the sum over $\al$ seems wasteful, the key observation is that if the partition of $[0,1]^{R+1}$ is chosen in an arithmetically meaningful way, this sum over $\al$ can later be precisely evaluated. This observation occurred first in \cite{HBP15} in the case of dimension $n=1$, and then  in \cite{Pie16} in arbitrary dimensions. The precise evaluation of the sum over $\al$, which we carry out in the next section, introduces bounds for the number of solutions to a system of Diophantine equations associated to $\Gbf$, and the corresponding analogue of the Vinogradov Mean Value Theorem. The known bound in the Vinogradov Mean Value Theorem leads to a savings that  compensates for the loss incurred by summing over all $\al$ in this step.

Finally, we remark on the fact that the right-hand side of (\ref{TS}) still contains a supremum, while we claimed our maneuvers aimed to remove the suprema from the objects we were considering. The point is that we will bound $S_F(\kbf)$ by a non-negative function that is increasing in the coordinates of $\kbf$, so that the supremum over $\kbf \leq 2\Hbf/P$ may be handled quite simply at a later step (see (\ref{SF_k})).

\section{Evaluation of the additive component}

We now turn to studying $S_F(\kbf)$ for a fixed $\kbf \leq 2\Hbf/P$. It is convenient to define the following notation. 
Given $2r$ tuples  $\xbf^{(1)}, \ldots, \xbf^{(2r)} \in \Z^{n}$, we will represent this collection by $\{\xbf\}$. For each $j=1,\ldots, 2r$, let $\ep(j) = (-1)^{j+1}$ and set $\del(j) = +1$ if $j$ is odd and $\Del-1$ if $j$ is even, where $\Del$ is the order of $\chi$ modulo $q$. Given such a collection $\{\xbf\}$, we then define
\[ \Sig_{\mathrm{add}} (\{ \xbf\}) := \sum_\al e \left( \sum_{j=1}^{2r} \ep(j) \theta_\al (\xbf^{(j)}) \right),\]
in which the sum over $\al$ denotes a sum over all the indices in the decomposition of $[0,1]^{R+1}$ constructed above.
Also define
\[ \Sig^F_{\mathrm{mult}}( \{\xbf\}) := \sum_{\mbf \modd{q}} \chi (F_{\{\xbf\}}(\mbf)),\]
in which 
\beq\label{dfn_Gpoly}
 F_{\{\xbf\}}(\Xbf) = \prod_{j=1}^{2r} F(\Xbf + \xbf^{(j)})^{\del(j)}.
 \eeq
Define $\Xi (\Gbf;\{\xbf\})$ to be the indicator function for  the set 
\[ V_r(\Gbf) :=  \{ \xbf^{(1)},\ldots,  \xbf^{(2r)} \in \mathbb{Z}^{n}  : \sum_{j=1}^{2r} \ep(j)( \xbf^{(j)})^\be = 0  , \forall  \be \in \Lambda(\Gbf) \}. \]
Later we will use the fact that  
\[ |V_r(\Gbf) \intersect (\zerobf,\kbf]^{2r} | \leq J_r(\Gbf,k_{\max}),\]
where  $k_{\max} = \max\{k_1,\ldots, k_n\}$,  
and we recall the notation that $J_r(\Gbf,X)$ counts the number of solutions of the system  (\ref{Vin_sys}).

Evaluating the sum over $\al$ in $\Sig_{\mathrm{add}} (\{ \xbf\})$ leads to the following identity, which we will apply with the choice $\Kbf = 2\Hbf/P$. Recall the weight $M=M(\Gbf)$ of the system $\Gbf$.
\begin{lemma}\label{lemma_prod}
Let $\Kbf = (K,\ldots, K)$. Upon setting $Q = \lceil 2r K \rceil,$   for each $\kbf \leq \Kbf$,
\beq\label{S_Sig}
  S_F(\kbf) =Q^M \sum_{\bstack{\xbf^{(1)},\ldots, \xbf^{(2r)} \in \Z^n}{ \mathbf{0} < \xbf^{(j)} \leq \kbf}}\Xi (\Gbf;\{\xbf\})\Sig^F_{\mathrm{mult}} (\{ \xbf\}).
  \eeq
\end{lemma}

To prove the lemma, expand the $2r$-th power in the definition of $S_F(\kbf)$ to show that 
\beq\label{S_AB}
 S_F(\kbf) =  \sum_{\bstack{\xbf^{(1)},\ldots, \xbf^{(2r)} \in \Z^n}{ \mathbf{0} < \xbf^{(j)} \leq \kbf}} \Sig_{\mathrm{add}}(\{\xbf\}) \Sig^F_{\mathrm{mult}} (\{ \xbf\}).
	\eeq
Now we recall  the partition $[0,1]^{R+1}	 = \union_\al B_\al$ and the distinguished vertices $\theta_\al$, which allow us to evaluate precisely the sum $\Sig_{\mathrm{add}}(\{\xbf\})$  for each fixed collection $\{\xbf\}$. 
Briefly (as also described  in \cite[\S 6]{Pie16}), by definition of the distinguished vertices $\theta_\al$ in (\ref{theta_spec}) and their associated polynomials,
\[ \Sig_{\mathrm{add}} (\{ \xbf\}) = \sum_\al e \left( \sum_{j=1}^{2r} \ep(j) \theta_\al (\xbf^{(j)}) \right)
	=\sum_{c_{\be^{(0)}}, \ldots, c_{\be^{(R)}}} e \left( \sum_{\be = \be^{(0)},\ldots, \be^{(R)}}
		c_\be \Qbf^{-\be} \left(   \sum_{j=1}^{2r} \ep(j) (\xbf^{(j)})^\be \right) \right),\]
		where the sum over $c_{\be^{(0)}}, \ldots, c_{\be^{(R)}}$ indicates summing for each $i =0,\ldots, R$ the parameter $c_{\be^{(i)}}$ over integers $0 \leq c_{\be^{(i)}} \leq \Qbf^{\be^{(i)}}-1$.
		This can be re-written as
\[ \Sig_{\mathrm{add}} (\{ \xbf\}) = 
	\prod_{\be =  \be^{(0)},\ldots, \be^{(R)}} 
		\left\{ \sum_{c_\be \modd{\Qbf^\be}} e \left( c_\be \Qbf^{-\be} \left(  \sum_{j=1}^{2r} \ep(j) (\xbf^{(j)})^\be \right) \right)
		\right\},\]
		so that by orthogonality of characters,
for each multi-index $\be$ we get a nonzero contribution of $\Qbf^\be = Q^{|\be|}$ if and only if  $\sum_{j=1}^{2r} \ep(j) (\xbf^{(j)})^\be \con 0$ modulo $\Qbf^\be$.
Precisely, we have shown
\[ \Sig_{\mathrm{add}}(\{\xbf\}) = Q^{M} \Xi_{\Qbf}(\Gbf;\{\xbf\}),\]
where we recall the definition of $M=M(\Gbf)$ from (\ref{rank_dfn}) and 
  we define $\Xi_{\Qbf}(\Gbf;\{\xbf\})$ to be the indicator function for  the set 
\beq\label{set}
  \{ \xbf^{(1)},\ldots,  \xbf^{(2r)} \in \mathbb{Z}^{n} : \sum_{j=1}^{2r} \ep(j)( \xbf^{(j)})^\be \con 0 \modd{\Qbf^\be}, \forall  \be \in \Lambda(\Gbf) \}. 
  \eeq
  \noindent  \emph{Remark:} Note that we do not need to consider a   congruence condition for $\Qbf^{\be}$ when $\be = (0,\ldots, 0)$, so we can write $\Lambda(\Gbf)$ instead of $\Lambda_0(\Gbf)$ in the definition of this set.
  
If  $Q$ is sufficiently large relative to $K$, we may force the congruences in the definition of this set (\ref{set}) to be identities in $\Z$, for every collection $\{\xbf\}$ such that $\xbf^{(j)} \in (\zerobf,\kbf]$, with $\kbf \leq \Kbf$. It suffices to choose 
\[Q  = \lceil 2rK \rceil .\]
Indeed, with this choice of $Q$, we see that for any fixed $\kbf \leq \Kbf$, each congruence in (\ref{set}) can only hold with $\xbf^{(j)} \in (\zerobf,\kbf]$  if it holds as an identity in $\Z$.
We conclude that 
\[ \Sig_{\mathrm{add}}(\{\xbf\}) = Q^M \Xi (\Gbf;\{\xbf\}),\]
and this proves Lemma \ref{lemma_prod}.

 \section{Stratification of the multiplicative component}
The next step is to count the number of collections $\{\xbf\}$ for which $\Sig^F_{\mathrm{mult}}(\{\xbf\})$ satisfies certain upper bounds. To do so, we will apply the stratification of Xu \cite{Xu18}, in the format of \cite[Theorem 4.4]{PieXu19}.
The key result we prove in this section is as follows; we will apply this with $\Kbf = 2\Hbf/P$.
 
\begin{prop}\label{prop_S}
Fix $n \geq 2$ and $r \geq n$. Suppose that $\Kbf = (K,\ldots, K)$ with $K \geq 1$ and $Q = \lceil 2rK \rceil$. 
Define $\Theta = \Theta_{n,r}= \lfloor (r-1)/(n-1) \rfloor$.
 Under the assumption $q^{1/2} K^{-\Theta} \leq 1$,  
 \[ \sup_{\kbf \leq \Kbf} S_F(\kbf) \ll  Q^M \{  J_r(\Gbf,K) q^{n/2}  +   K^{2nr-\Theta} q^{n/2+1/2}  \}.\]
 \end{prop}
 \noindent
\emph{Remark:} In the standard terminology of the Burgess method, the first term in braces may be seen as the contribution of ``good'' collections $\{ \xbf\}$, namely those that lead to a complete character sum $\Sig^F_{\mathrm{mult}}(\{\xbf\})$ with square-root cancellation. A key result of Xu's work is that such ``good'' collections are generic among all tuples in $\Z^{2nr} \intersect (\zerobf,\kbf]^{2r}$, once $r$ is sufficiently large that $\Theta_{n,r} \geq 1$. This term includes a factor $ J_r(\Gbf,K)$ instead of $K^{2nr}$ because of the advantageous evaluation of the additive character sum, leading to the presence of the indicator function $\Xi(\Gbf; \{\xbf\})$ in (\ref{S_Sig}), which imposes that the collection $\{\xbf\}$ must lie in the set $V_r(\Gbf)$. The savings of $J_r(\Gbf,K)$ relative to $K^{2nr}$ will compensate for the large factor $Q^M$ in front, which we accrued by summing over all indices $\al$ during the proof of Lemma \ref{lemma_TS}.
 
The second term in braces is the contribution of the ``bad'' collections $\{\xbf\}$. 
The  ``bad'' collections lead to character sums with bounds ranging from $O(q^{(n+1)/2})$ to $O(q^n)$;   Xu's stratification helpfully shows that  ``bad'' collections have positive codimension in $\Z^{2nr}$, and the collections that yield progressively worse bounds for the character sum have progressively higher codimension.

 Before we prove Proposition \ref{prop_S}, let us see how it implies Proposition \ref{prop_T}.
 We   apply the bound from Proposition \ref{prop_S} with the choice $K=2H/P$ in (\ref{TS}),   so that $Q = \lceil 4rH/P \rceil$ and so $Q^M \ll_r (H/P)^M$. 
 We conclude that 
\[
 T(F,\Gbf; \Nbf,\Hbf) \ll  (H/P)^{M/2r} H^{-n/2r}P^{n-1/2r}q^{n/4r}( \log q)^{n+1} 	 
 \{  J_r(\Gbf,2H/P)^{1/2r}  + q^{1/4r} (H/P)^{n- \Theta/2r}   \},
\]
   as claimed in Proposition \ref{prop_T}.

 \subsection{Proof of Proposition \ref{prop_S}}

 We first define some notation. 
 Fix $n \geq 2, r \geq n$. Recall that $\Theta = \Theta_{n,  r} = \lfloor \frac{r-1}{n-1} \rfloor.$ For any   $1 \leq j \leq n$, and for a tuple $\kbf = (k_1,\ldots, k_n)$ with $k_1 \leq \cdots \leq k_n$, we define 
\[
 B_{n,r}(j;\kbf) =  B_{n,r}(j;k_1, \ldots, k_n) 
 	=\begin{cases}
 1 & j=0 \\
  k_1^{j\Theta} =  k_1^{j\lfloor \frac{r-1}{n-1} \rfloor}  & j=1,\ldots, n-2 \\
 k_1^{r-1}  & j=n-1 \\
  (k_1 \cdots k_{n/2})^{2r}  & j=n,  \text{$n$ even} \\
    (k_1 \cdots k_{(n-1)/2})^{2r}k_{(n+1)/2}^{r}  & j=n, \text{$n$ odd} . 
 \end{cases}
\]
  We now recall \cite[Thm. 4.4]{PieXu19}, which is essentially the result of \cite{Xu18}, specialized to our setting.
 \begin{letterthm}\label{thm_Xu}
Let integers $n \geq 2$, $r\geq n$, $\Del \geq 2$, $D \geq 1$ be fixed.  Then there exist constants $C=C(n,r,D) \geq 1$ and $C'' = C''(n,r,\Del,D) \geq 1$ such that the following holds.

Fix a prime $q$, and let $\chi$ be a non-principal multiplicative Dirichlet character of conductor $q$ and order $\Del$. Let $F \in \Z[x_1,\ldots, x_n]$  be a form of degree $D$ with  $(\Del,q)$-admissible reduction modulo $q$. Define $F_{\{\xbf\}} (\Xbf)$ for each collection $\{\xbf\} \in \Z^{2nr}$ as in (\ref{dfn_Gpoly}).
Then for every  $1\le j\le n$, for every tuple $\kbf = (k_1,\ldots, k_n) \in \Z^n$ with $1 \leq k_1 \leq k_2 \leq \cdots \leq k_n \leq q$, 
\beq\label{Xu_ineq}
 \# \left\{ (\xbf^{(1)}, \ldots, \xbf^{(2r)})
\in (\zerobf,\kbf]^{2r} : \left|  \sum_{\mbf \modd{q}}
\chi(F_{\{\xbf\}}(\mbf))
\right|  > C q^{(n+j-1)/2} \right\}
 \leq C''  \| \kbf \|^{2r} B_{n,r}(j;\kbf)^{-1}.
 \eeq
\end{letterthm}

\noindent  \emph{Remark:} 
For more information on the stratification theorem of Xu in this setting, see \cite[\S2, \S4]{PieXu19}, as well as the original work \cite{Xu18}. Roughly speaking, the exponent $\Theta = \Theta_{n, r}$ arises from a lower bound on the codimension of a subscheme of those collections $\{\xbf\} = (\xbf^{(1)},\ldots, \xbf^{(2r)}) \in \Z^{2nr}$ for which square-root cancellation could fail to hold for the complete character sum.
The number of collections with a corresponding complete character sum that exceeds square-root cancellation, i.e. $C q^{n/2}$, is bounded above by $C''  \| \kbf \|^{2r} k_1^{-\Theta}$, the number of collections with a corresponding complete character sum  that exceeds $C q^{(n+1)/2}$ is bounded above by the smaller quantity $C''  \| \kbf \|^{2r} k_1^{-2\Theta}$, and so on. This motivates the definition of the functions  $B_{n, r}(j; \kbf)$.

Now we prepare to apply the stratification result of Theorem \ref{thm_Xu} to $S_F(\kbf)$. Let $C=C(n,r,D)$ be the constant provided by Theorem \ref{thm_Xu}.
Let us fix $\kbf \leq \Kbf,$ not yet assuming the ordering $k_1 \leq \ldots \leq k_n $.
 For each $1 \leq j \leq n$, define 
\beq\label{Yj}
Y_j= Y_j^F(\kbf):= \left\{ \{\xbf\} 
\in (\zerobf,\kbf]^{2r} : \left|  \sum_{\mbf \modd{q}}
\chi(F_{\{\xbf\}}(\mbf))
\right|  > C q^{(n+j-1)/2} \right\}.
\eeq
Then $(\zerobf,\kbf]^{2r} =:Y_0\supset Y_1\supset Y_2\supset\dots\supset Y_n\supset Y_{n+1}:=\varnothing$.
 Upon employing the disjoint dissection $(\zerobf,\kbf]^{2r}=\coprod_{j=0}^n Y_j\setminus Y_{j+1}$ in (\ref{S_Sig}),
we now see that for this fixed $\kbf$,
\begin{align*}
 S_F(\kbf) &  \leq   Q^M \sum_{j=0}^n \sum_{\{ \xbf\} \in (Y_j \setminus Y_{j+1}) \intersect V_r(\Gbf) } |\Sig^F_{\mathrm{mult}} (\{ \xbf\})|\\
 	& \leq Q^M \sum_{j=0}^n  \# (Y_j \intersect V_r(\Gbf)) C q^{(n+(j+1)-1)/2 }.
	\end{align*}
Given that the method leading to Theorem \ref{thm_Xu} (see \cite{Xu18})  can only compute upper bounds for $\# Y_j$ in terms of the dimension of $Y_j$, it is difficult to obtain a nontrivial upper bound for the intersection $Y_j \intersect V_r(\Gbf)$, except in the case of $Y_0 = (\zerobf,\kbf]^{2r}$.
But in this case, we see that  
\[  \# (Y_0 \intersect V_r(\Gbf)) = \# (V_r(\Gbf)\intersect (\zerobf,\kbf]^{2r} ) \leq J_r(\Gbf,k_{\max}).\]
Thus we obtain 
\beq\label{SF_k_0}
 S_F(\kbf) 
\leq C Q^M  J_r(\Gbf,k_{\max}) q^{n/2}  + C Q^M \sum_{j=1}^n  \# (Y_j^F(\kbf))  q^{(n+(j+1)-1)/2 }.
\eeq
At this point, if in particular $k_1 \leq k_2 \leq \cdots k_n \leq  q$ then we can apply
	Theorem \ref{thm_Xu} in the form of the upper bound 
 $\#Y_j^F(\kbf) \leq C'' \|\kbf\|^{2r}B_{n,r}(j;\kbf)^{-1}$ for each $1 \leq j \leq n$, uniformly in $F$.
Consequently in this case we have
\beq\label{SF_k}
   S_F(\kbf) \leq C Q^M J_r(\Gbf,k_n) q^{n/2}  + C  C'' Q^M  \sum_{j=1}^n  q^{(n+j)/2} \|\kbf\|^{2r}B_{n,r}(j;\kbf)^{-1} .
   \eeq
 More generally, given any fixed $\kbf$, we will re-order the variables $x_1,\ldots, x_n$  in $F(x_1,\ldots, x_n)$ so that the correct ordering does hold for the entries in $\kbf$, 
 and then we will apply Theorem \ref{thm_Xu} to a form defined according to this re-ordering. This will use the uniformity of the bound in Theorem \ref{thm_Xu}, with respect to the form $F$.
 
We now give the precise argument. Given any permutation $\pi$ on $\{1,\ldots, n\}$, define the form $F^{\pi}(x_1,\ldots, x_n)  = F(x_{\pi(1)}, x_{\pi(2)}, \ldots, x_{\pi(n)})$; $F^\pi$ has  $(\Del,q)$-admissible reduction modulo $q$  if and only if $F$ does.  Given any fixed $\kbf \leq \Kbf$, 	  let $\sig$ be a permutation on the indices $\{1,\ldots, n\}$ such that 
\beq\label{k_order}
 k_{\sig(1)} \leq k_{\sig(2)} \leq \cdots \leq k_{\sig(n)}.\eeq
Given any tuple $\xbf \in \Z^{n}$ let $\xbf_\sig = (x_{\sig(1)},\ldots, x_{\sig(n)})$ and similarly let $\{\xbf\}_\sig$ denote the collection $\xbf^{(1)}_\sig, \ldots, \xbf^{(2r)}_\sig.$
Since $F(x_1,\ldots, x_n) = F^{\sig^{-1}} (x_{\sig(1)}, \ldots, x_{\sig(n)})$ we see that for any fixed collection $\{\xbf\}$,
\[   \sum_{\mbf \modd{q}} \chi(F_{\{\xbf\}}(\mbf)) 
	=  \sum_{\mbf_\sig \modd{q}} \chi(F^{\sig^{-1}}_{\{\xbf\}_\sig}(\mbf_\sig))
	= \sum_{\mbf \modd{q}} \chi(F^{\sig^{-1}}_{\{\xbf\}_\sig}(\mbf)).
\]
Of course, the set of all $\{ \xbf\}_\sig \in (\zerobf,\kbf_\sig]^{2r} $ identifies with
the set of all $\{ \xbf\}  \in (\zerobf,\kbf]^{2r} $.
In particular, if we recall that $Y_j^F(\kbf)$ denotes the set in (\ref{Yj}) and let 
\beq\label{Yj'}
Y_j^{F^{\sig^{-1}}}(\kbf_\sig):= \left\{ \{\xbf\}_{\sig} 
\in (\zerobf,\kbf_\sig]^{2r} : \left|  \sum_{\mbf \modd{q}}
\chi(F^{\sig^{-1}}_{\{\xbf\}_\sig}(\mbf))
\right|  > C q^{(n+j-1)/2} \right\},
\eeq
 we see that 
$\#Y_j^F(\kbf) = \#Y_j^{F^{\sig^{-1}}}(\kbf_\sig)$ for each $j=1,\ldots, n$. 
We will apply this inside each term in (\ref{SF_k_0}) (which we recall holds for $\kbf$ without assuming an ordering on the coordinates of $\kbf$). We thus obtain from (\ref{SF_k_0}) that 
\[ 
 S_F(\kbf) 
\leq  C Q^M  J_r(\Gbf,k_{\max}) q^{n/2}  + CQ^M \sum_{j=1}^n  \# (Y_j^{F^{\sig^{-1}}}(\kbf_\sig) )  q^{(n+(j+1)-1)/2 }.
\]
 Now $\kbf_\sig$ satisfies the ordering (\ref{k_order}) and thus we may apply Theorem \ref{thm_Xu} to bound the cardinality of the sets $Y_j^{F^{\sig^{-1}}}(\kbf_\sig)$ and  conclude that 
\beq\label{SF_k}
   S_F(\kbf) \leq C Q^M J_r(\Gbf,k_{\max}) q^{n/2}  + C  C''   Q^M\sum_{j=1}^n  q^{(n+j)/2} \|\kbf_\sig\|^{2r}B_{n,r}(j;\kbf_\sig)^{-1} .
   \eeq
 
 Observe that with respect to the variable $\kbf \in \R^n$, $ \|\kbf\|^{2r}B_{n,r}(j;\kbf)^{-1}$ is a non-decreasing function in each coordinate of $\kbf$; that is, for fixed $r \geq 1$, for each $1 \leq j \leq n$, there exist  exponents $\al_{j,1}, \ldots, \al_{j,n} \geq 0$ (also depending on $r$) such that 
 \[ 
 \|\kbf\|^{2r}B_{n,r}(j;\kbf)^{-1}= k_1^{\al_{j,1}} \cdots k_n^{\al_{j,n}},
 \]
for all tuples $\kbf$ with $k_1 \leq \cdots \leq k_n$. This is an immediate consequence of the definition of the functions $B_{n, r}(j; \cdot)$.
In particular, for any  $\kbf \leq \Kbf$ with $k_1 \leq \cdots \leq k_n$, where $\Kbf = (K, \ldots, K)$, we obtain that 
\[ 
\|\kbf\|^{2r}B_{n,r}(j;\kbf)^{-1} \leq K^{\al_{j,1}} \cdots K^{\al_{j,n}}= \|\Kbf\|^{2r}B_{n,r}(j;\Kbf)^{-1} .
\]
We apply this to each term in (\ref{SF_k}), and we conclude that 
\[ \sup_{\kbf \leq \Kbf} S_F(\kbf) \ll_{n, r, D, \Del}   Q^M J_r(\Gbf,K) q^{n/2}  + Q^M  \|\Kbf\|^{2r} q^{n/2}\sum_{j=1}^n  q^{j/2} B_{n,r}(j;\Kbf)^{-1}.\]
 Proposition \ref{prop_S}  follows, after we verify a  lemma about sums of the functions $B_{n, r}(j; \cdot)$,
 which we prove in a general setting.
\begin{lemma}\label{lemma_B_sum}
If $1 \leq K_1 \leq K_2 \leq \cdots \leq K_n$ then for each $r \geq n$,
\beq\label{B_sum}
\sum_{j=1}^n  q^{j/2} B_{n,r}(j;\Kbf)^{-1}\ll q^{1/2} K_1^{-\Theta} 
\eeq
as long as \beq\label{K_relation}
q^{1/2} K_1^{-\Theta} \leq 1.
\eeq
\end{lemma}
 
By  definition, the sum over $j=1,\ldots, n-1$ takes the form
\beq\label{sum_j}
     \sum_{j=1}^{n-2}q^{j/2} K_1^{-j\Theta}  + q^{(n-1)/2}K_1^{-(r-1)}    \leq   \sum_{j=1}^{n-1}(q^{1/2} K_1^{-\Theta}  )^j
\eeq
in which $\Theta = \lfloor (r-1)/(n-1) \rfloor$; this used the fact that  $r-1 \geq (n-1) \lfloor (r-1)/(n-1) \rfloor = (n-1)\Theta$. The right-most expression shows that under the assumption (\ref{K_relation}), all terms $j \geq 2$ are dominated by $j=1$. 
The last term to check is $j=n$. For $n \geq 2$ even, the $j=n$ term is 
\[q^{n/2}(K_1 \cdots K_{n/2})^{-2r}  \leq q^{n/2} K_1^{-nr} \leq (q^{1/2}K_1^{-\Theta})^n, \]
in which we have used the ordering $K_1 \leq K_2 \leq \cdots \leq K_n$ and the fact that $r \geq (r-1)/(n-1) \geq \Theta$  when $n \geq 2$. 
For $n \geq 3$ odd, similar reasoning shows the $j=n$ term is 
	\[q^{n/2}(K_1 \cdots K_{(n-1)/2})^{-2r}K_{(n+1)/2}^{-r} \leq q^{n/2} K_1^{-rn}\leq (q^{1/2}K_1^{-\Theta})^n.  \]
 In either case, under the assumption (\ref{K_relation}) we see that the $j=n$ term is $\ll q^{1/2} K_1^{-\Theta}$.  
This completes the proof of the lemma, verifying Proposition   \ref{prop_S}, and hence Proposition \ref{prop_T}.

\noindent  \emph{Remark:} On the other hand, observe that the  summands on the right-hand side of (\ref{sum_j}) increase with $j$ if $q^{1/2} K_1^{-\Theta} >1$; this will motivate our later choice of $P$ so that (\ref{K_relation}) holds, with $K_1= 2H/P$.

     \section{Concluding arguments for Theorem \ref{thm_mixed}}\label{sec_conc}
     With Proposition \ref{prop_T} in hand, the final steps to prove Theorem \ref{thm_mixed} are to apply a bound for $J_r(\Gbf,X)$ and to choose $P$.
    In order to motivate our choice for $P \leq H$, we recall that so far we have supposed in (\ref{HP_assp}) and the application of Proposition \ref{prop_S} with $K=2H/P$ that
      \beq\label{P_cond}
HP<q  , \qquad \qquad   P \leq H q^{-1/2\Theta}.
      \eeq
   
      We first argue formally in some generality, in order to understand the role of the Vinogradov Mean Value Theorem.
      We suppose we are in a range of $r$ where 
      \beq\label{Jr_assp}
      J_r(\Gbf,X) \ll X^{2rn-\mu+\ep},
      \eeq
       for some positive integer $\mu$ (depending on $n,d, r$ and the system $\Gbf$). As we remark below,
       this is   known for $\mu =M = M(\Gbf)$ for all values of $r$ that we will consider, but for later reference we initially argue  in terms of the abstract parameter $\mu$.

        Under these assumptions,
\[ J_r(\Gbf,2H/P)^{1/2r}  +  (H/P)^{n-\Theta/2r} q^{1/4r}  
\ll q^\ep((H/P)^{n-\mu/2r}   + (H/P)^{n- \Theta/2r} q^{1/4r} ).\]
	Now we observe that to balance these two terms we would choose $P$ to be an integer with 
\beq\label{P_choose}
 \frac{1}{2} Hq^{-\frac{1}{2(\Theta-\mu)}} \leq    P < Hq^{ -\frac{1}{2(\Theta-\mu)}} .
 \eeq
This supposes that $\Theta > \mu$ in order to meet the requirement that $P \leq H$; since $\Theta = \lfloor (r-1)/(n-1) \rfloor$, this is a requirement that $r$ is sufficiently large with respect to $\mu,n$. 
This choice for $P$ also satisfies the requirements   in (\ref{P_cond}), as long as we assume that $H<q^{\frac{1}{2}+\frac{1}{4 (\Theta - \mu)}}$. 

\noindent  \emph{Remark:} This will be satisfied, by a hypothesis of the theorem, when we ultimately apply this reasoning with $\mu = M = M(\Gbf)$.

We apply this choice of $P$ in Proposition \ref{prop_T} to conclude that if (\ref{Jr_assp}) holds then
\beq\label{T_penultimate}
  T(F,\Gbf; \Nbf,\Hbf) \ll H^{n-\frac{n+1}{2r}} q^{\frac{n(\Theta - \mu) + M+1 - \mu}{4r(\Theta-\mu)}}q^\ep,
\eeq
for any $\ep>0$,
as long as $ \Theta = \lfloor (r-1)/(n-1) \rfloor  > \mu.$ This is a condition on $r$, namely $r > \mu + 1$ when $n=2$, and in general it suffices to have 
$r > (\mu +1)(n-1) + 1$ when $n \geq 3$. 

Now to understand $\mu$ in (\ref{Jr_assp}), we  restrict our attention to $\Gbf$ being the standard system (\ref{G_choice}) for dimension $n$ and degree $d$. (We remark on more general systems in \S \ref{sec_ranges} in an appendix.)
One can calculate that in order for $\mu, r$ to be such that (\ref{Jr_assp}) holds and simultaneously $ \Theta = \lfloor (r-1)/(n-1) \rfloor  > \mu,$
we must be in the range of $r$ such that the savings in (\ref{Jr_assp}) is   $\mu = M(\Gbf)$. 
(We provide the details to prove this simple observation in the appendix.)
 This comes from the known upper bounds in the multi-dimensional Vinogradov Mean Value Theorem, which we now recall.
 
Precisely, for all $n \geq 2$, (\ref{Jr_assp}) with $\mu = M$ is true for all values of $r$ satisfying
\beq\label{rM}
r > M+1 \quad (n=2), \qquad  r > (M+1)(n-1) + 1 \quad (n \geq 3),
\eeq
 due to the truth of the Vinogradov Mean Value Theorem for the system $\Gbf$ defined in (\ref{G_choice}).
To be precise, for $n\geq 3$, $J_{r}(\Gbf,X)\ll X^{2nr - M+ \ep}$ for $r$ in the range (\ref{rM}) is known from \cite{PPW13}.
\noindent  \emph{Remark:} This uses the fact that for $n \geq 3$, the requirement on $r$ in (\ref{rM}) imposes that $r > R(d+1)$, which was the requirement in the work \cite{PPW13}.

On the other hand, for $n=2$, for each $d \geq 2$, in order to obtain this upper bound for $r$ in the range (\ref{rM}) one requires the stronger results of \cite{GuoZha18}, which apply for all $r \geq 1$.

Thus we   now  only consider  the case that (\ref{Jr_assp}) holds with $\mu = M=M(\Gbf)$ and $r$ is in the range (\ref{rM}).
When both these conditions are met,   (\ref{T_penultimate}) shows that
\beq\label{T_penultimate'}
  T(F,\Gbf; \Nbf,\Hbf) \ll H^{n-\frac{n+1}{2r}} q^{\frac{n(\Theta - M) +1  }{4r(\Theta- M)}}q^\ep.
\eeq
 This suffices to complete the proof of Theorem \ref{thm_mixed},  since $|S(F, g;\Nbf,\Hbf)| \leq  T(F,\Gbf; \Nbf,\Hbf)$.

\subsection{Quantification of the strength of Theorem \ref{thm_mixed}}\label{sec_quant}
Supposing that $H=q^\be$, then the bound for $|S(F, g \; \Nbf,\Hbf)|$  provided by Theorem \ref{thm_mixed} 
is nontrivial i.e. $o(q^\be)$ as long as 
\beq\label{beta_thres}
\be > \frac{1}{2} - \frac{(\Theta - M-1)}{2(\Theta -M)(n+1)},
\eeq
in which $\Theta = \lfloor \frac{r-1}{n-1} \rfloor.$
This allows for values of $\be$ strictly smaller than  $1/2$ as long as $r$ is sufficiently large that $\Theta > M+1$. The right-hand side in (\ref{beta_thres}) is always $ > \be_n$ with $\be_n$ as defined in (\ref{beta_dfn}). We thus suppose that $H = q^{\be_n + \kappa}$ for some small $\kappa>0$, in which case we can compute that 
$|S(F, g;\Nbf,\Hbf)| \ll H^n q^{-\del}$ with 
\[ \del = \frac{2\kappa (n+1) (\Theta - M) -1}{4r (\Theta - M)}.
\]
We now use the approximation of replacing $\Theta$ by $(r-1)/(n-1)$ (which in fact is exact, when $n=2$, and not far off from the truth when $r$ grows very large, as it will when we choose $r$ according to $\kappa$ and $\kappa \maps 0$).
After this approximation, we can write $\del$ as the value at $r$ of the function 
\[ f_{b,c,d}(r) = \frac{ br-c}{r(r-d)}\]
with 
\[ 
b= \frac{\kappa (n+1)}{2}, \qquad c= (M(n-1)+1) \frac{\kappa (n+1)}{2} + \frac{n-1}{4}, \qquad 
	d = M(n-1)+1.\]
The function $f_{b,c,d}(r)$ attains a local extremum  at $r=b^{-1}(c \pm \sqrt{c^2 - dbc})$; using the values for 
$b,c,d$ above and simplifying using $\kappa \maps 0$ we see that we should choose $r$ to be the nearest integer to 
\[ r \approx \frac{n-1}{n+1} \cdot \frac{1}{\kappa}.\]
This choice of $r$ satisfies $\Theta = \Theta_{n, r} > M$ if $\kappa$ is sufficiently small (relative to $n,d$). 
We now apply this in the expression above for $\del$, now further approximating $\Theta$ by $r/(n-1)$, and we see that in the limit as $\kappa \maps 0$ we obtain a savings over $H^n$ of the form  $H^n q^{-\del}$, in which
\[ \del \approx \frac{(n+1)^2}{4(n-1)} \kappa^2.\]
The significance of this savings is that it is independent of the degree $d$ of $g$, due to  the application of the multi-dimensional Vinogradov Mean Value Theorem. In particular, it is as strong as the savings of the first author and Xu \cite{PieXu19} in the purely multiplicative case.

\section{Appendix: Further remarks on Vinogradov systems}\label{sec_gen}
In this appendix, we briefly remark on three aspects of the proof of Theorem \ref{thm_mixed}. First, we consider how improvements to   $\Theta$   would lead to a setting in which one would require the sharp results of \cite{GuoZha18} for all $n \geq 3$, in addition to the results of \cite{PPW13}. Second, we explain why we only considered $\mu=M$ in the conclusion of the proof of Theorem \ref{thm_mixed}, or equivalently, why the current Burgess method in this setting leads to consideration of very large $r$. Third, we briefly state a more general result for systems $\Gbf$ other than the standard system (\ref{G_choice}).

\subsection{Remarks on the codimension $\Theta$}
Let $\Gbf$ denote the standard system (\ref{G_choice}) of monomials in $n$ variables of degree at most $d$. In the proof of Theorem \ref{thm_mixed}, we applied the results of \cite{PPW13} to bound $J_r(\Gbf,X)$ when $n \geq 3$, and only required the stronger results of \cite{GuoZha18}  when $n=2$. We now remark that if one could improve the stratification of Xu for complete multiplicative character sums (in the sense of the discussion in  \cite[\S 8.2]{PieXu19}), then one would require the results of \cite{GuoZha18} for all $n \geq 2$.
 
Precisely, we have seen in the   argument in \S \ref{sec_conc} that if the known upper bound is $J_r(\Gbf,X) \ll X^{2nr-\mu + \ep}$, then the result of the  Burgess method developed in this paper must restrict to values of $r$ for which $\Theta > \mu$. Here $\Theta$ is the codimension of the first exceptional subscheme  $X_1$ arising in the stratification of Xu \cite{Xu18}; see the remark following Theorem \ref{thm_Xu} for a rough idea. Currently, for $n \geq 2$, Xu has obtained  $\Theta \geq \lfloor (r-1)/(n-1) \rfloor$. 
The appendix in \cite[\S 8.2]{PieXu19} outlines  conjectural possibilities for  improvements to the codimension leading to the value for $\Theta$;  for example, one might hope to prove that $\Theta=r$ is possible.
 In dimension $n=2$ this is nearly attained already by $\Theta = \lfloor (r-1)/(n-1) \rfloor$, but is significantly different from the current result for large $n$. Let us suppose that one could prove $\Theta  \geq r/\al(n)$ for some function $1 \leq \al(n) \leq n-1$, leading to the restriction $r > \mu \al(n)$ in the method of proof for Theorem \ref{thm_mixed}. In particular, if we use $\mu = M$ and $\al(n)$ is not too large, then in order to obtain our theorem unconditionally for all $r> M \al(n)$ we would require  the results of \cite{GuoZha18} for $J_r(\Gbf,X)$ for those $r$ with $M \al (n) < r < R(d+1)$, while  \cite{PPW13} would continue to apply for $r> R(d+1)$.

\subsection{Remarks on intermediate ranges of $r$}\label{sec_ranges}
In the proof of Theorem \ref{thm_mixed}, we remarked that we need only consider  the upper bound (\ref{Jr_assp})  when $\mu = M$, 
  where $M = M(\Gbf)$ is the weight of the associated system of Diophantine equations (the sum of the total degrees), and $r$ is very large. One might ask whether one could consider other values for $\mu$, and correspondingly smaller values $r$. 
  Here we explain why the Burgess method developed in this paper only allows the regime of $r$ in which $\mu =M =M(\Gbf)$. 
  
  In the current discussion we can take $\Gbf$ to be any reduced monomial translation-dilation invariant system.
  Precisely, the question is: what must $\mu$ be in order for both (\ref{Jr_assp}) and $\Theta= \lfloor (r-1)/(n-1) \rfloor  > \mu$ to hold?
Given any reduced monomial translation-invariant system $\Gbf$, suppose there is a sequence of positive integers $K_j = K_j(\Gbf)$ for $1 \leq j \leq n$  such that for all $X \geq 1$, for all $r \geq 1$, 
\beq\label{J_bound_upper}
J_{r}(\Gbf,X) \ll_{r,\Gbf,\ep} X^\ep (  X^{rn} + \sum_{j=1}^{n} X^{2rj+(n-j) - K_{j}}).
\eeq
(Note that in this notation, $K_n$ plays the role of $M(\Gbf)$.)
In particular, by the breakthrough work of Guo and Zhang \cite{GuoZha18}, this is now known for the standard system $\Gbf$ in (\ref{G_choice}) of monomials in  $n$ variables with total degree at most $d$; in this case
$K_j=  \frac{jd}{j+1} {j+d \choose j}$. 

 We claim that if a bound of the form (\ref{J_bound_upper}) holds, then in order for both (\ref{Jr_assp}) and $\lfloor (r-1)/(n-1) \rfloor > \mu$ to hold (and hence certainly $r > \mu$), we must have $\mu = K_n$ (so that the $j=n$ term   dominates in (\ref{J_bound_upper})). 
Indeed, suppose that $r$ is such that the $j$-th term dominates in (\ref{J_bound_upper}),  for some $1 \leq j \leq n$. In the notation of (\ref{Jr_assp}), this would impose $\mu = (2r-1)(n-j) + K_j$. Then in order to have  $r > \mu$ we must at least have $r > (2r-1)(n-j)$, which   can only hold if $n=j$. 
(Similarly, the term $X^{rn}$ cannot dominate,  since that would impose $\mu = rn$, but the condition $r > \mu$ could not hold.) This proves the claim.  (Even if, for example, the codimension $\Theta$ could be improved to $r$, the analogue of (\ref{P_choose}) would still require $r> \mu$, leading to $\mu = M(\Gbf)$ via the same argument given above.)
 Thus it appears that significant innovations to the method would be required, in order to be able to apply counts for Vinogradov systems where any term with $j< n$ dominates in (\ref{J_bound_upper}).
 
\subsection{Remarks on other systems} 
Let $\Gbf$ be any reduced monomial translation-dilation invariant system, in any dimension $n \geq 2$ and with degree $d(\Gbf) \geq 1$.
  Parsell, Prendiville and Wooley proved that (\ref{J_bound_upper}) holds for any $r > R(\Gbf) (d(\Gbf)+1)$, in which case the $j=n$ term dominates (with $K_n(\Gbf)=M(\Gbf)$), and the upper bound is 
$J_r(\Gbf,X) \ll X^{2rn-M(\Gbf)+ \ep}$. 
More recently, Guo and Zorin-Kranich \cite{GuoZor18} have proved that a sharp upper bound of the form (\ref{J_bound_upper}), with appropriately defined $K_j (\Gbf)$, holds for all $r \geq 1$, for more general systems $\Gbf$, which we now describe. 
Fix a tuple $(k_1,\ldots, k_n)$ of positive integers. Fix an integer $k$. Let $\Gbf$ be the system defined 
according to the set of exponents   
\beq\label{G_ACK}
\Lam(\Gbf )= \{ \be: \be_1 \leq k_1, \ldots, \be_n \leq k_n , 1 \leq |\be| \leq k\}.
\eeq
If in particular $k=k_1 = \cdots = k_n = d$ then this is the standard system (\ref{G_choice}).
If $k=k_1 + \cdots + k_n$ then this is known as an Arkhipov-Chubarikov-Karatsuba system.
 For the systems Guo and Zorin-Kranich handle we can thus   obtain a generalization of Theorem \ref{thm_mixed}   in the largest range of $r$ allowed by the Burgess method developed in this paper.
 We record the conclusion of this discussion:
\begin{thm}\label{thm_mixed_T}
Fix $n \geq 2$ and $d,D \geq 1$. Let $q$ be a fixed prime, and let $\chi$ be a non-principal Dirichlet character of order $\Del$ and conductor $q$. 
Let $F \in \Z[x_1,\ldots, x_n]$ be a form of degree $D$ such that its reduction modulo $q$ is $(\Del,q)$-admissible. 
Let $\Gbf$ be a reduced monomial translation-dilation invariant system with rank $R(\Gbf)$ and weight $M(\Gbf)$ (containing linear monomials in each variable).
For each integer $r \geq 1$, define $\Theta = \lfloor (r-1)/(n-1) \rfloor$. 
Then
\beq\label{thm_ineq}
  |T(F,\Gbf;\Nbf,\Hbf) | \ll  H^{n-\frac{n+1}{2r}} q^{\frac{n(\Theta - M(\Gbf)) +1 }{4r(\Theta-M(\Gbf))}}q^\ep, 
  \eeq
  for every integer $r$ such that $r > R(\Gbf)(d(\Gbf)+1)$ and $\Theta > M(\Gbf)$,
  and for every $\Hbf = (H, \ldots, H)$ with $H< q^{1/2 + 1/(4(\Theta - M(\Gbf)))}$. Furthermore if $\Gbf$ is a system of the type (\ref{G_ACK}) then we may take any $r$ such that $\Theta> M(\Gbf)$. 
The implied constant could depend on $\Gbf, n,D, \Del, r,\ep$ but is otherwise independent of $F$.
\end{thm}

\section*{Acknowledgements}
Pierce is partially supported by NSF   CAREER grant DMS-1652173, a Sloan Research Fellowship, and the AMS Joan and Joseph Birman Fellowship.  Pierce thanks the Hausdorff Center for Mathematics  and the initiative ``A Room of One's Own'' for focused environments, and the referee for helpful comments. 
Pierce also thanks the nine individuals who helped to tend her children during the main time period in which this work took place.

%***************************************
\bibliographystyle{alpha}
\bibliography{NoThBibliography}
%***************************************
\label{endofproposal}

\end{document}